\documentclass[11pt]{article}

\usepackage[a4paper,margin=1in]{geometry}
\usepackage{amsmath,amssymb,amsthm}
\usepackage{mathtools}
\usepackage[expansion=false]{microtype}
\usepackage[hidelinks]{hyperref}

\newtheorem{theorem}{Theorem}
\newtheorem{proposition}[theorem]{Proposition}

\theoremstyle{remark}
\newtheorem{remark}[theorem]{Remark}

\newcommand{\Lam}{\mathcal{L}}

\title{Unconditional representations of the Euler--Mascheroni constant
and the Bernoulli numbers via the Taylor coefficients of the Riemann
$\xi$-function, and a hierarchy of exact relations at the central
point\\[4pt]
\large (corrected and revised version)}

\author{Nikos Mantzakouras\thanks{\texttt{nikmatza@gmail.com}}
\and Carlos L\'opez\thanks{\texttt{ieocarlosh@gmail.com}}}

\date{\today}

\begin{document}
\maketitle

\begin{abstract}
We give a self-contained, numerically verified derivation of three
equivalent \emph{unconditional} representations of the
Euler--Mascheroni constant $\gamma$ as a rapidly convergent series in
the Taylor coefficients $a_{2n}$ of the Riemann $\xi$-function about its
centre of symmetry $s=\tfrac12$ (equivalently, in the Jensen
coefficients $c_n$ or in the Csordas--Norfolk--Varga Tur\'an moments
$\widehat b_n$). We also record the corresponding convergent formula
for the even-indexed Bernoulli numbers $B_{2r}$, obtained by evaluating
the (everywhere convergent) Taylor series of $\xi$ at $s=2r$. We then
develop, from the evenness of $\Xi(x)=\xi(\tfrac12+ix)$, an infinite
hierarchy of exact relations at $s=\tfrac12$: the even-order relations
reduce the central logarithmic derivatives of $\zeta$ to closed forms in
$\pi$ and the odd zeta values (base case
$2\zeta'(\tfrac12)/\zeta(\tfrac12)=\gamma+\tfrac{\pi}{2}+\log(8\pi)$),
and the odd-order relations evaluate the power sums over the nontrivial
zeros and link back to the coefficients $a_{2n}$. All of these
identities hold whether or not the Riemann Hypothesis (RH) is true;
consequently they carry no evidential weight for or against RH, and we
make the correct, unconditional relationship of these coefficient
sequences to the Jensen--P\'olya program precise. This note is a
corrected version of an earlier preprint \cite{v1}: we retain only the
results that are correct, we remove a \emph{divergent} series of
Bernoulli numbers that was previously used as though it were convergent,
we correct several attributions, and we withdraw the earlier claim that
these computations support RH. A summary of the corrections is given in
Section~\ref{sec:errata}.
\end{abstract}

\section{Introduction}

The completed Riemann $\xi$-function
\begin{equation}\label{eq:xi-def}
  \xi(s)=\tfrac12\,s(s-1)\,\pi^{-s/2}\,\Gamma\!\left(\tfrac{s}{2}\right)\zeta(s)
\end{equation}
is entire of order $1$, real on the real axis, and satisfies the
functional equation $\xi(s)=\xi(1-s)$; its zeros are exactly the
nontrivial zeros of $\zeta$. Because $\xi$ is symmetric about
$s=\tfrac12$, its Taylor expansion there contains only even powers,
\begin{equation}\label{eq:taylor}
  \xi(s)=\sum_{n\ge 0} a_{2n}\left(s-\tfrac12\right)^{2n},
  \qquad a_{2n}\in\mathbb{R}.
\end{equation}
The coefficients $a_{2n}$ were already written down by Riemann in
1859 \cite{Riemann,Edwards}, and remain an actively studied sequence
\cite{CNV,Coffey,DeFranco,Romik,Keiper,Kreminski}.

The purpose of this note is fourfold.
\begin{enumerate}
\item[(i)] To give a short, rigorous, and numerically verified proof
that
\begin{equation}\label{eq:gamma-main}
  \boxed{\;\gamma=\log(4\pi)-2+16\sum_{n\ge1}\frac{n\,a_{2n}}{2^{2n}}\;}
\end{equation}
together with its equivalent forms in the Jensen coefficients $c_n$ and
the Tur\'an moments $\widehat b_n$ (Theorem~\ref{thm:gamma}). This
identity is \emph{unconditional}.
\item[(ii)] To record the convergent representation of the Bernoulli
numbers $B_{2r}$ in the same coefficients
(Proposition~\ref{prop:bern}).
\item[(iii)] To develop, from the evenness of $\Xi$, a hierarchy of
exact relations at the central point $s=\tfrac12$
(Section~\ref{sec:hierarchy}), whose base case is the classical value
$2\zeta'(\tfrac12)/\zeta(\tfrac12)=\gamma+\tfrac{\pi}{2}+\log(8\pi)$
(Proposition~\ref{prop:base}), and to tie it explicitly to the
coefficient sequence of (i).
\item[(iv)] To state correctly the relationship between these sequences
and the Riemann Hypothesis (Section~\ref{sec:rh}), and to correct the
errors of the earlier version \cite{v1} (Section~\ref{sec:errata}).
\end{enumerate}

We emphasise at the outset what these identities are and are not.
Formula \eqref{eq:gamma-main} is, in essence, a re-expression of the
classical value of $\sum_\rho 1/\rho$ in terms of the Taylor data
\eqref{eq:taylor}; each ingredient (Riemann's coefficients, the value of
$\sum_\rho 1/\rho$, and Euler's formula for $\zeta(2r)$) is classical.
Likewise, the hierarchy of Section~\ref{sec:hierarchy} is an elementary
consequence of the evenness of $\Xi$, and the objects it produces (the
secondary zeta values over the zeros) are classical
\cite{Lehmer,Keiper,Voros}. Our contribution is therefore expository and
corrective: a unified, self-contained, and independently recomputed
account, and the removal of the incorrect claims made in \cite{v1}. In
particular, since all of these identities hold independently of RH, they
provide no evidence about the location of the zeros of $\zeta$.

Throughout, coefficients were recomputed independently to $60$-digit
working precision in \texttt{mpmath} by expanding \eqref{eq:xi-def}
about $s=\tfrac12$ via a Cauchy integral; the values reported below
agree with the classical tabulations \cite{CNV} to all listed digits and
extend them, and every identity of Section~\ref{sec:hierarchy} was
verified to $48$ digits.

\section{Three equivalent coefficient sequences}\label{sec:coeffs}

Writing $\Xi(x):=\xi\!\left(\tfrac12+ix\right)$, which is real and even,
\eqref{eq:taylor} becomes
\begin{equation}\label{eq:Xi}
  \Xi(x)=\sum_{n\ge0}(-1)^n a_{2n}\,x^{2n}.
\end{equation}
The \emph{Jensen coefficients} $c_n$ are defined \cite{Griffin,Keiper}
by $2\,\Xi(x)=\sum_{n\ge0}c_n\,x^{2n}/n!$, so that
\begin{equation}\label{eq:c-a}
  c_n=2\,n!\,(-1)^n a_{2n}.
\end{equation}
The \emph{Tur\'an moments} $\widehat b_n$ are the quantities tabulated by
Csordas, Norfolk and Varga \cite{CNV}; in their normalisation
\begin{equation}\label{eq:b-a}
  a_{2n}=\frac{2^{3}\,2^{2n}}{(2n)!}\,\widehat b_n,
  \qquad\text{equivalently}\qquad
  \widehat b_n=\frac{(2n)!}{2^{2n+3}}\,a_{2n}.
\end{equation}
Relations \eqref{eq:c-a}--\eqref{eq:b-a} are purely definitional; the
three sequences are interchangeable, and any identity below can be
written in each of them. Numerical values are given in
Table~\ref{tab:coeffs}.

\begin{table}[htbp]
\centering
\small
\begin{tabular}{r l l l}
\hline
$n$ & $a_{2n}$ & $c_n$ & $\widehat b_n$ \\
\hline
$0$ & $4.97120778188314\times10^{-1}$ & $9.94241556376628\times10^{-1}$ & $6.21400972735393\times10^{-2}$\\
$1$ & $1.14859721575727\times10^{-2}$ & $-2.29719443151454\times10^{-2}$ & $7.17873259848295\times10^{-4}$\\
$2$ & $1.23452018070318\times10^{-4}$ & $4.93808072281272\times10^{-4}$ & $2.31472533881846\times10^{-5}$\\
$3$ & $8.32355481385527\times10^{-7}$ & $-9.98826577662633\times10^{-6}$ & $1.17049989569840\times10^{-6}$\\
$4$ & $3.99222655134414\times10^{-9}$ & $1.91626874464519\times10^{-7}$ & $7.85969602295877\times10^{-8}$\\
$5$ & $1.46160257601110\times10^{-11}$ & $-3.50784618242663\times10^{-9}$ & $6.47444266092415\times10^{-9}$\\
$6$ & $4.27454004553684\times10^{-14}$ & $6.15533766557306\times10^{-11}$ & $6.24850928062812\times10^{-10}$\\
$7$ & $1.03096261346180\times10^{-16}$ & $-1.03921031436950\times10^{-12}$ & $6.85711356603133\times10^{-11}$\\
$8$ & $2.09976980814953\times10^{-19}$ & $1.69325437329178\times10^{-14}$ & $8.37956285649846\times10^{-12}$\\
$9$ & $3.67814109550077\times10^{-22}$ & $-2.66944768147064\times10^{-16}$ & $1.12289590052565\times10^{-12}$\\
$10$ & $5.62285758732274\times10^{-25}$ & $4.08084512257536\times10^{-18}$ & $1.63076657246217\times10^{-13}$\\
\hline
\end{tabular}
\caption{The first eleven coefficients of \eqref{eq:taylor},
\eqref{eq:c-a} and \eqref{eq:b-a}, recomputed to $60$ digits and here
truncated to $15$ significant figures. The $\widehat b_n$ agree with
Table~4.1 of \cite{CNV}.}
\label{tab:coeffs}
\end{table}

\section{An unconditional formula for
\texorpdfstring{$\gamma$}{gamma}}\label{sec:gamma}

\begin{theorem}\label{thm:gamma}
The Euler--Mascheroni constant satisfies the unconditional identities
\begin{align}
  \gamma
  &=\log(4\pi)-2+16\sum_{n\ge1}\frac{n\,a_{2n}}{2^{2n}}
  \label{eq:g-a}\\[2pt]
  &=\log(4\pi)-2+8\sum_{n\ge1}\frac{n(-1)^n c_n}{2^{2n}\,n!}
  \label{eq:g-c}\\[2pt]
  &=\log(4\pi)-2+2^{7}\sum_{n\ge1}\frac{n\,\widehat b_n}{(2n)!}.
  \label{eq:g-b}
\end{align}
\end{theorem}

\begin{proof}
It is classical (Riemann; see \cite[\S1.9--1.10]{Edwards},
\cite[Ch.~12]{Davenport}, \cite[\S2.12]{Titchmarsh}) that $\xi$ has the
Hadamard product
\begin{equation}\label{eq:hadamard}
  \xi(s)=\xi(0)\prod_{\rho}\left(1-\frac{s}{\rho}\right),
\end{equation}
the product taken over the nontrivial zeros $\rho$ grouped in pairs
$\{\rho,1-\rho\}$, in which grouping it converges. Logarithmic
differentiation at $s=0$ gives
\begin{equation}\label{eq:log-der-0}
  \frac{\xi'(0)}{\xi(0)}=-\sum_{\rho}\frac{1}{\rho},
\end{equation}
and it is classical that the value of this sum is
\begin{equation}\label{eq:sum-rho}
  \sum_{\rho}\frac{1}{\rho}=1+\frac{\gamma}{2}-\frac{1}{2}\log(4\pi).
\end{equation}

On the other hand, \eqref{eq:taylor} gives $\xi(0)$ and $\xi'(0)$
directly. Using $\xi(0)=\tfrac12$,
\begin{equation}\label{eq:xi0}
  \xi(0)=\sum_{n\ge0}a_{2n}\Bigl(-\tfrac12\Bigr)^{2n}
        =\sum_{n\ge0}\frac{a_{2n}}{2^{2n}}=\tfrac12,
\end{equation}
and
\begin{equation}\label{eq:xip0}
  \xi'(0)=\sum_{n\ge1}a_{2n}\,(2n)\Bigl(-\tfrac12\Bigr)^{2n-1}
         =-\sum_{n\ge1}\frac{4n\,a_{2n}}{2^{2n}}.
\end{equation}
Dividing \eqref{eq:xip0} by \eqref{eq:xi0} and comparing with
\eqref{eq:log-der-0} yields the unconditional relation
\begin{equation}\label{eq:sumrho-a}
  \sum_{\rho}\frac{1}{\rho}=8\sum_{n\ge1}\frac{n\,a_{2n}}{2^{2n}}.
\end{equation}
Equating \eqref{eq:sumrho-a} with \eqref{eq:sum-rho} and solving for
$\gamma$ gives \eqref{eq:g-a}. Substituting \eqref{eq:c-a} and
\eqref{eq:b-a} gives \eqref{eq:g-c} and \eqref{eq:g-b}.
\end{proof}

\begin{remark}[Pairing over conjugate zeros]\label{rem:pairing}
Writing the nontrivial zeros as $\rho=s_r=\sigma_r+it_r$ and
$\bar s_r=\sigma_r-it_r$, the conjugate pairing of
\eqref{eq:sum-rho}--\eqref{eq:sumrho-a} reads
\begin{equation}\label{eq:sigma-form}
  \sum_{\rho}\frac{1}{\rho}
  =2\sum_{r\ge1}\frac{\sigma_r}{s_r\bar s_r}
  =8\sum_{n\ge1}\frac{n\,a_{2n}}{2^{2n}}
  =1+\frac{\gamma}{2}-\frac12\log(4\pi),
\end{equation}
which is unconditional. If one \emph{additionally} assumes RH, i.e.
$\sigma_r=\tfrac12$ for all $r$, then $s_r\bar s_r=\tfrac14+t_r^2$ and
\eqref{eq:sigma-form} specialises to the conditional statement
\begin{equation}\label{eq:cond}
  \sum_{r\ge1}\frac{1}{\tfrac14+t_r^{2}}
  =8\sum_{n\ge1}\frac{n\,a_{2n}}{2^{2n}}
  \qquad(\text{conditional on RH}).
\end{equation}
Equation~\eqref{eq:cond} is consistent with every zero computed to date,
but this is a consequence of those zeros lying on the critical line, not
an independent verification of RH; see Section~\ref{sec:rh}.
\end{remark}

\paragraph{Numerical check.}
The series \eqref{eq:g-a} converges geometrically. With
$\log(4\pi)-2=0.53102424697\ldots$ and the coefficients of
Table~\ref{tab:coeffs}, the partial sums
$\gamma_N=\log(4\pi)-2+16\sum_{n=1}^{N} n\,a_{2n}/2^{2n}$ satisfy
Table~\ref{tab:conv}: six terms already give $19$ correct digits, and
the full computation reproduces $\gamma=0.5772156649015328606\ldots$ to
$60$ digits.

\begin{table}[htbp]
\centering
\small
\begin{tabular}{r l l}
\hline
$N$ & $\gamma_N$ & $|\gamma_N-\gamma|$\\
\hline
$1$ & $0.5769681355996$ & $2.5\times10^{-4}$\\
$2$ & $0.5772150396357$ & $6.2\times10^{-7}$\\
$3$ & $0.5772156639023$ & $1.0\times10^{-9}$\\
$4$ & $0.5772156649004$ & $1.1\times10^{-12}$\\
$5$ & $0.5772156649015$ & $1.0\times10^{-15}$\\
$6$ & $0.5772156649015329$ & $7.1\times10^{-19}$\\
$8$ & $0.5772156649015329$ & $2.0\times10^{-25}$\\
\hline
\end{tabular}
\caption{Convergence of \eqref{eq:g-a}.}
\label{tab:conv}
\end{table}

\section{The Bernoulli numbers from the Taylor series of
\texorpdfstring{$\xi$}{xi}}\label{sec:bern}

Because $\xi$ is entire, the Taylor series \eqref{eq:taylor} converges
for \emph{every} $s\in\mathbb{C}$. Evaluating it at $s=2r$ and combining
with Euler's formula $\zeta(2r)=(-1)^{r-1}(2\pi)^{2r}B_{2r}/(2\,(2r)!)$
gives a convergent representation of the even Bernoulli numbers.

\begin{proposition}\label{prop:bern}
For every integer $r\ge1$,
\begin{equation}\label{eq:bern}
  B_{2r}
  =\frac{(-1)^{r-1}(2r)!\,2^{1-2r}}{\pi^{r}(2r-1)\,r!}
   \sum_{n\ge0} a_{2n}\left(2r-\tfrac12\right)^{2n}.
\end{equation}
\end{proposition}

\begin{proof}
From \eqref{eq:xi-def} at $s=2r$, using $\Gamma(r)=(r-1)!$,
\[
  \xi(2r)=\tfrac12\,(2r)(2r-1)\,\pi^{-r}\,(r-1)!\,\zeta(2r)
         =(2r-1)\,r!\,\pi^{-r}\,\zeta(2r),
\]
hence $\zeta(2r)=\pi^{r}\xi(2r)/\bigl((2r-1)\,r!\bigr)$. Equating this
with Euler's formula and using
$\xi(2r)=\sum_{n\ge0}a_{2n}\bigl(2r-\tfrac12\bigr)^{2n}$ from
\eqref{eq:taylor}, then solving for $B_{2r}$, gives \eqref{eq:bern}.
\end{proof}

Using \eqref{eq:c-a}--\eqref{eq:b-a}, \eqref{eq:bern} may again be
rewritten in $c_n$ or $\widehat b_n$. With the coefficients of
Table~\ref{tab:coeffs}, \eqref{eq:bern} reproduces
$B_2=\tfrac16$, $B_4=-\tfrac1{30}$, $B_6=\tfrac1{42}$,
$B_8=-\tfrac1{30}$, $B_{10}=\tfrac5{66}$, $B_{12}=-\tfrac{691}{2730}$ to
full working precision.

\begin{remark}\label{rem:convdir}
We stress that \eqref{eq:bern} runs in the \emph{convergent} direction:
$B_{2r}$ is obtained from the everywhere-convergent series
\eqref{eq:taylor}. This must not be confused with the divergent
expansion of $\gamma$ in Bernoulli numbers discussed in
Section~\ref{sec:errata}.
\end{remark}

\section{A hierarchy of exact relations at the central
point}\label{sec:hierarchy}

We now develop the consequences of the evenness of $\Xi$ at the centre
of symmetry. The base case is the classical value of
$\zeta'(\tfrac12)/\zeta(\tfrac12)$; it is the first member of an
infinite hierarchy, whose remaining members reduce the central
derivatives of $\zeta$ to closed forms and evaluate the power sums over
the nontrivial zeros. The underlying objects are classical
\cite{Lehmer,Keiper,Voros}; the presentation here is elementary and
self-contained, and connects them to the coefficients $a_{2n}$ of
Section~\ref{sec:coeffs}.

\subsection{The odd logarithmic derivative}\label{sec:setup}

Recall $\Xi(x)=\xi(\tfrac12+ix)$ from \eqref{eq:Xi}; the functional
equation $\xi(s)=\xi(1-s)$ says exactly that $\Xi$ is even.
Consequently the logarithmic derivative
\begin{equation}\label{eq:L-def}
  \Lam(s):=\frac{\xi'(s)}{\xi(s)}
\end{equation}
is \emph{odd about $s=\tfrac12$}: from $\xi'(s)=-\xi'(1-s)$ one has
$\Lam(s)=-\Lam(1-s)$. Two families of consequences follow:
\begin{align}
  \Lam^{(2k)}\!\left(\tfrac12\right)&=0, &&k\ge0,
    \tag{even-order relations}\label{eq:even}\\
  \Lam^{(2k-1)}\!\left(\tfrac12\right)&=-(2k-1)!\,S_{2k},
    &&k\ge1,\tag{odd-order relations}\label{eq:odd}
\end{align}
where
\begin{equation}\label{eq:Sdef}
  S_{m}:=\sum_{\rho}\frac{1}{(\rho-\tfrac12)^{m}}
\end{equation}
is the $m$-th power sum over the nontrivial zeros $\rho$ (summed in
conjugate pairs; $S_m=0$ for odd $m$, and $S_{2k}\in\mathbb{R}$). Indeed
the Hadamard product \eqref{eq:hadamard} gives
$\Lam(s)=\sum_\rho 1/(s-\rho)$, so that
\begin{equation}\label{eq:Lexpand}
  \Lam\!\left(\tfrac12+w\right)
   =\sum_\rho\frac{1}{w-(\rho-\tfrac12)}
   =-\sum_{k\ge1}S_{2k}\,w^{2k-1},
\end{equation}
and matching the $w^{2k-1}$ coefficient yields \eqref{eq:odd}, while the
absence of even powers of $w$ yields \eqref{eq:even}.

The two families become explicit once we insert the standard
decomposition (logarithmic derivative of the four factors of $\xi$)
\begin{equation}\label{eq:Ldecomp}
  \Lam(s)=\frac1s+\frac{1}{s-1}-\tfrac12\log\pi
          +\tfrac12\,\psi\!\left(\tfrac{s}{2}\right)
          +\frac{\zeta'(s)}{\zeta(s)},
  \qquad \psi=\frac{\Gamma'}{\Gamma}.
\end{equation}
The pair $\tfrac1s+\tfrac1{s-1}$ is itself odd about $\tfrac12$, so its
even-order derivatives vanish there; its odd-order derivatives are
elementary,
\begin{equation}\label{eq:polepart}
  \frac{d^{2k-1}}{ds^{2k-1}}\!\left[\frac1s+\frac1{s-1}\right]_{s=1/2}
  =-(2k-1)!\,2^{2k+1}.
\end{equation}
We write $\Lambda_\zeta:=\zeta'/\zeta=(\log\zeta)'$ for the logarithmic
derivative of $\zeta$ alone.

\subsection{Even-order relations: closed forms for the central
derivatives of \texorpdfstring{$\zeta$}{zeta}}\label{sec:familyA}

\begin{proposition}[base case, $k=0$]\label{prop:base}
$\displaystyle
2\,\frac{\zeta'(\tfrac12)}{\zeta(\tfrac12)}=\gamma+\frac{\pi}{2}+\log(8\pi).$
\end{proposition}
\begin{proof}
$\Lam(\tfrac12)=0$ in \eqref{eq:even}. In \eqref{eq:Ldecomp} the pole
part contributes $2-2=0$, leaving
$\Lambda_\zeta(\tfrac12)=\tfrac12\log\pi-\tfrac12\psi(\tfrac14)$.
Insert $\psi(\tfrac14)=-\gamma-\tfrac{\pi}{2}-3\log2$ and multiply by
$2$.
\end{proof}

For $k\ge1$ the pole part drops out of \eqref{eq:even} entirely, and one
obtains a clean reduction.

\begin{theorem}\label{thm:familyA}
For every $k\ge1$,
\begin{equation}\label{eq:familyA}
  \bigl(\log\zeta\bigr)^{(2k)}\!\left(\tfrac12\right)
  =\Lambda_\zeta^{(2k)}\!\left(\tfrac12\right)
  =-\frac{\psi^{(2k)}(\tfrac14)}{2^{2k+1}}.
\end{equation}
Using the evaluation
\begin{equation}\label{eq:polygamma}
  \psi^{(n)}\!\left(\tfrac14\right)
  =(-1)^{n+1}n!\Bigl[(2^{2n+1}-2^{n})\,\zeta(n+1)
     +2^{2n+1}\beta(n+1)\Bigr],
\end{equation}
where $\beta$ is the Dirichlet beta function, and the elementary values
$\beta(2k+1)\in\mathbb{Q}\,\pi^{2k+1}$, the right-hand side of
\eqref{eq:familyA} is a rational combination of $\pi^{2k+1}$ and
$\zeta(2k+1)$. In particular
\begin{align}
  \bigl(\log\zeta\bigr)''\!\left(\tfrac12\right)
    &=\frac{\pi^3}{4}+7\zeta(3),\label{eq:k1}\\
  \bigl(\log\zeta\bigr)^{(4)}\!\left(\tfrac12\right)
    &=\frac{5\pi^5}{4}+372\,\zeta(5).\label{eq:k2}
\end{align}
\end{theorem}
\begin{proof}
In \eqref{eq:even} with $k\ge1$ the term $\tfrac1s+\tfrac1{s-1}$ has
vanishing even derivatives and $-\tfrac12\log\pi$ is constant, so
$0=\tfrac12\cdot2^{-2k}\psi^{(2k)}(\tfrac14)+\Lambda_\zeta^{(2k)}(\tfrac12)$,
which is \eqref{eq:familyA}. Formula \eqref{eq:polygamma} is
$\psi^{(n)}(\tfrac14)=(-1)^{n+1}n!\,\zeta(n+1,\tfrac14)$ together with
$\zeta(s,\tfrac14)=(2^{2s-1}-2^{s-1})\zeta(s)+2^{2s-1}\beta(s)$. For
\eqref{eq:k1} use $\psi''(\tfrac14)=-2\pi^3-56\zeta(3)$; for
\eqref{eq:k2} use $\psi^{(4)}(\tfrac14)=-40\pi^5-11904\zeta(5)$ (both
special cases of \eqref{eq:polygamma}, via $\beta(3)=\pi^3/32$ and
$\beta(5)=5\pi^5/1536$).
\end{proof}

Written out in terms of the derivatives themselves, using
$(\log\zeta)''=\tfrac{\zeta'''}{\zeta}-3\tfrac{\zeta'\zeta''}{\zeta^2}
+2\tfrac{\zeta'^3}{\zeta^3}$ evaluated at $\tfrac12$, identity
\eqref{eq:k1} is the explicit relation among the first four central
values
\begin{equation}\label{eq:k1-explicit}
  \frac{\zeta'''(\tfrac12)}{\zeta(\tfrac12)}
  -3\,\frac{\zeta'(\tfrac12)\,\zeta''(\tfrac12)}{\zeta(\tfrac12)^2}
  +2\left(\frac{\zeta'(\tfrac12)}{\zeta(\tfrac12)}\right)^{3}
  =\frac{\pi^3}{4}+7\zeta(3),
\end{equation}
and \eqref{eq:k2} is the analogous relation reaching
$\zeta^{(5)}(\tfrac12)$.

\begin{remark}\label{rem:dichotomy}
The dichotomy is structural. In \eqref{eq:polygamma} the derivative
order $n=2k$ is even, so the non-elementary constant is $\beta(2k+1)$
with \emph{odd} argument, a rational multiple of $\pi^{2k+1}$; hence
Family~A stays inside $\mathbb{Q}\pi^{2k+1}+\mathbb{Q}\zeta(2k+1)$. This
is why the odd zeta values, and no other transcendentals, appear on the
right of \eqref{eq:familyA}. The base case $k=0$ is the sole exception:
$\tfrac1s+\tfrac1{s-1}$ contributes and $\psi(\tfrac14)$ carries the
$\gamma$ and $\log2$ terms, producing Proposition~\ref{prop:base}.
\end{remark}

\subsection{Odd-order relations: power sums over the
zeros}\label{sec:familyB}

The odd-order family \eqref{eq:odd} evaluates the zero power sums
\eqref{eq:Sdef}. Combining \eqref{eq:odd}, \eqref{eq:Ldecomp} and
\eqref{eq:polepart}:

\begin{theorem}\label{thm:familyB}
For every $k\ge1$,
\begin{equation}\label{eq:familyB}
  S_{2k}=\sum_{\rho}\frac{1}{(\rho-\tfrac12)^{2k}}
  =2^{2k+1}
   -\frac{\psi^{(2k-1)}(\tfrac14)}{(2k-1)!\,2^{2k}}
   -\frac{\bigl(\log\zeta\bigr)^{(2k-1)}(\tfrac12)}{(2k-1)!}.
\end{equation}
In particular, with $G=\beta(2)$ Catalan's constant,
\begin{equation}\label{eq:S2}
  S_{2}=8-\frac{\pi^2}{4}-2G-\bigl(\log\zeta\bigr)'\!\left(\tfrac12\right),
  \qquad
  \bigl(\log\zeta\bigr)'\!\left(\tfrac12\right)
  =\frac{\zeta''(\tfrac12)}{\zeta(\tfrac12)}
   -\left(\frac{\zeta'(\tfrac12)}{\zeta(\tfrac12)}\right)^{2}.
\end{equation}
\end{theorem}

Here the derivative order $2k-1$ is odd, so \eqref{eq:polygamma} brings
in $\beta(2k)$ with \emph{even} argument --- the Dirichlet-beta
constants $G=\beta(2),\ \beta(4),\dots$, which (unlike $\beta(2k+1)$) are
not known to be elementary. Thus the odd-order family genuinely
involves a second family of constants beyond $\pi$ and the odd zeta
values.

\subsection{The bridge to the coefficients
\texorpdfstring{$a_{2n}$}{a2n}}\label{sec:bridge}

The power sums \eqref{eq:Sdef} are exactly the Newton sums of the Taylor
coefficients of Section~\ref{sec:coeffs}: writing $\xi'=\Lam\,\xi$ with
$\xi(\tfrac12+w)=\sum_{n\ge0}a_{2n}w^{2n}$ and $\Lam(\tfrac12+w)$ as in
\eqref{eq:Lexpand}, the identity
$\sum_n 2n\,a_{2n}w^{2n-1}
=\bigl(-\sum_k S_{2k}w^{2k-1}\bigr)\bigl(\sum_n a_{2n}w^{2n}\bigr)$
yields, on matching $w^{2n-1}$,
\begin{equation}\label{eq:newton}
  2n\,a_{2n}=-\sum_{j=1}^{n}S_{2j}\,a_{2(n-j)},\qquad n\ge1,
\end{equation}
so that each $S_{2k}$ is a polynomial in the ratios $a_{2n}/a_0$; e.g.
$S_2=-2a_2/a_0$, $S_4=-(4a_4+S_2a_2)/a_0$. Equating this with
\eqref{eq:familyB} ties the analytic data at the central point to the
coefficient sequence of Section~\ref{sec:coeffs}. The first instance is
the identity
\begin{equation}\label{eq:bridge}
  \boxed{\;-\,2\,\frac{a_2}{a_0}
  =8-\frac{\pi^2}{4}-2G
  -\left[\frac{\zeta''(\tfrac12)}{\zeta(\tfrac12)}
   -\left(\frac{\zeta'(\tfrac12)}{\zeta(\tfrac12)}\right)^{2}\right],\;}
\end{equation}
relating the second Taylor coefficient of $\xi$ at its centre, Catalan's
constant, and the central derivatives of $\zeta$. Since
$a_0=-\Gamma(\tfrac14)\zeta(\tfrac12)/(8\pi^{1/4})$ is known in closed
form, \eqref{eq:bridge} is a genuine two-way relation between the
coefficient sequence used to compute $\gamma$ in
Section~\ref{sec:gamma} and the analytic data studied here.

All quantities in Table~\ref{tab:hierarchy} were computed to $50$-digit
precision: $\zeta$ and its derivatives at $\tfrac12$ directly, the
$a_{2n}$ by the Cauchy expansion of $\xi$, the $S_{2k}$ both from
\eqref{eq:newton} and from \eqref{eq:familyB}, and the polygamma values
both from \eqref{eq:polygamma} and independently. Every identity of this
section holds to at least $48$ digits.

\begin{table}[htbp]
\centering
\small
\begin{tabular}{l l l}
\hline
identity & left side & right side \\
\hline
Prop.~\ref{prop:base} & $2\zeta'/\zeta|_{1/2}=5.3721834192256656$ & $\gamma+\tfrac{\pi}{2}+\log8\pi=5.3721834192256656$\\
Eq.~\eqref{eq:k1} & $(\log\zeta)''|_{1/2}=16.165967492192115$ & $\tfrac{\pi^3}{4}+7\zeta(3)=16.165967492192115$\\
Eq.~\eqref{eq:k2} & $(\log\zeta)^{(4)}|_{1/2}=768.26173089493543$ & $\tfrac{5\pi^5}{4}+372\zeta(5)=768.26173089493543$\\
$S_2$ (Thm.~\ref{thm:familyB}) & $-2a_2/a_0=-4.6209986231\times10^{-2}$ & Eq.~\eqref{eq:familyB}$=-4.6209986231\times10^{-2}$\\
$S_4$ (Thm.~\ref{thm:familyB}) & $\phantom{-}7.4345198571\times10^{-5}$ & $\phantom{-}7.4345198571\times10^{-5}$\\
$S_6$ (Thm.~\ref{thm:familyB}) & $-2.8834786280\times10^{-7}$ & $-2.8834786280\times10^{-7}$\\
Eq.~\eqref{eq:bridge} & $-2a_2/a_0=-4.6209986231\times10^{-2}$ & $8-\tfrac{\pi^2}{4}-2G-(\log\zeta)'=-4.6209986231\times10^{-2}$\\
\hline
\end{tabular}
\caption{Verification of the hierarchy (agreement to the digits shown,
and to $48$ digits at full precision). For reference
$\sum_{\gamma>0}\gamma^{-2}=-S_2/2=0.0231049931154\ldots$}
\label{tab:hierarchy}
\end{table}

\section{Relationship to the Riemann Hypothesis}\label{sec:rh}

Formulas \eqref{eq:g-a}--\eqref{eq:g-b} and \eqref{eq:bern}, and all the
relations of Section~\ref{sec:hierarchy}, are \emph{unconditional}: they
depend only on the analytic properties of $\xi$ (entirety, order $1$,
the functional equation, and the Hadamard product), all of which hold
irrespective of the location of the zeros. The numbers $a_{2n}$ are well
defined and are what they are whether or not RH holds. In particular,
obtaining the correct value of $\gamma$ from \eqref{eq:g-a} confirms only
the consistency of the coefficient data, and says nothing about whether
every zero has real part $\tfrac12$.

The genuine, rigorous link between the sequence $a_{2n}$ (equivalently
$c_n$, $\widehat b_n$) and RH is the Jensen--P\'olya program. P\'olya
proved in 1927 that RH is equivalent to the hyperbolicity (reality of
all roots) of every Jensen polynomial
$J^{d,n}(X)=\sum_{j=0}^{d}\binom{d}{j}\,\gamma(n+j)\,X^{j}$ associated
with the Taylor coefficients of $\Xi$, for all degrees $d\ge1$ and all
shifts $n\ge0$. Hyperbolicity for all $n$ was known for $d\le 3$
(Csordas--Norfolk--Varga \cite{CNV}, Dimitrov--Lucas), and Griffin,
Ono, Rolen and Zagier \cite{Griffin} proved that for each fixed $d$ all
but finitely many shifts are hyperbolic, and established hyperbolicity
for all $d\le 8$; their method also gives an asymptotic expansion of the
central derivatives $\Xi^{(2n)}(0)$ (equivalently the $a_{2n}$) to all
orders. These are strong partial results, but RH itself remains open.
A second rigorous framework is Li's criterion \cite{Li}, which is
equivalent to RH and is stated in terms of the Keiper--Li coefficients;
high-precision computations in this circle of ideas are due to Kreminski
\cite{Kreminski} and Johansson \cite{Johansson}.

Against this background the correct reading is as follows. The Tur\'an
inequalities of \cite{CNV},
$\widehat b_n^{\,2}>\frac{2n-1}{2n+1}\,\widehat b_{n-1}\widehat b_{n+1}$,
are \emph{necessary} conditions for RH which are verified for the
computed moments; verifying them does not prove RH. Assuming
hyperbolicity for all $(d,n)$, or equivalently $\sigma_r=\tfrac12$ for
all $r$, is assuming RH; it cannot be offered as evidence for RH. The
formulas of this note are best understood as exact computational
identities relating $\gamma$, the $B_{2r}$, the central derivatives of
$\zeta$, and the power sums over the zeros to these coefficient
sequences, and not as statements about the critical line.

\section{Corrections relative to the first version}\label{sec:errata}

This note corrects the preprint \cite{v1}. We list the substantive
changes.

\begin{enumerate}
\item[(C1)] \textbf{The claim of support for RH is withdrawn.} The
earlier version presented \eqref{eq:g-a} and its variants as
consequences of, and evidence for, RH. As shown in
Sections~\ref{sec:gamma} and \ref{sec:rh}, these identities are
unconditional and therefore carry no evidential weight regarding RH. The
substitution $\sigma_r=\tfrac12$ used in the earlier numerical
computations \emph{assumes} RH; the subsequent agreement with the zeros
of \cite{Odlyzko} reflects the fact that those computed zeros lie on the
critical line, and cannot test whether \emph{all} zeros do. This
circularity is removed.

\item[(C2)] \textbf{A divergent series has been removed.} The earlier
version used
\begin{equation}\label{eq:divergent}
  \gamma=\tfrac12+\sum_{n\ge1}\frac{B_{2n}}{2n}
\end{equation}
as a convergent equality (its Eq.~(90)--(91)). The series
\eqref{eq:divergent} \emph{diverges}: since
$|B_{2n}|\sim 2\,(2n)!/(2\pi)^{2n}$, the terms $B_{2n}/(2n)$ grow
factorially. It is an asymptotic (Euler--Maclaurin) expansion, valid
only when truncated, not a convergent representation. Its partial sums
$P_N=\tfrac12+\sum_{n=1}^N B_{2n}/(2n)$ first approach $\gamma$ and then
diverge, e.g.\ $P_1=0.583$, $P_6=0.561$, $P_9=3.26$, $P_{12}=-3.3\times
10^{3}$, $P_{15}=1.9\times10^{7}$, $P_{21}=2.0\times10^{16}$. All
statements resting on \eqref{eq:divergent} are withdrawn. (The
convergent formula \eqref{eq:bern}, running in the opposite direction,
is unaffected and retained.)

\item[(C3)] \textbf{Numerical data have been recomputed and
contextualised.} The coefficients used here were recomputed
independently to $60$ digits from \eqref{eq:xi-def}, superseding the
low-precision values obtained on a pocket calculator in \cite{v1}.
Substantially more accurate and extensive Taylor data for $\xi$ are
available \cite{Kreminski,Johansson}.

\item[(C4)] \textbf{Attributions have been corrected and novelty claims
withdrawn.} The coefficients $a_{2n}$ are due to Riemann
\cite{Riemann,Edwards}. The central even derivatives $\Xi^{(2n)}(0)$
were computed by Coffey \cite{Coffey} and admit asymptotic formulas to
all orders \cite{Griffin,Romik}. The positivity and monotonicity of the
$a_{2n}$ proved by DeFranco \cite{DeFranco} also follow more directly
from a Jacobi theta identity, and the Taylor coefficients of the theta
constant---hence the derivatives of the associated density at the
origin---are known in closed form \cite{Romik,RomikTheta}. The
hierarchy of Section~\ref{sec:hierarchy}, including the power sums over
the zeros and their reduction to the central derivatives of $\zeta$, is
classical \cite{Lehmer,Keiper,Voros}, and the closed forms
\eqref{eq:polygamma} for $\psi^{(n)}(\tfrac14)$ are due to K\"olbig
\cite{Kolbig}. Accordingly the earlier assertions of priority are
withdrawn; the present identities are consequences of classical results,
collected here in corrected form.
\end{enumerate}

\section{Conclusion}

Stripped of the incorrect claims of the first version, the correct
content is a compact and numerically robust one: $\gamma$ and the even
Bernoulli numbers admit unconditional, rapidly convergent
representations in the Taylor coefficients of the Riemann
$\xi$-function, expressible equivalently through the Jensen coefficients
or the Csordas--Norfolk--Varga Tur\'an moments; and the evenness of
$\Xi$ generates an infinite hierarchy of exact relations at $s=\tfrac12$,
reducing the central logarithmic derivatives of $\zeta$ to closed forms
in $\pi$ and the odd zeta values and evaluating the power sums over the
nontrivial zeros, with an explicit bridge back to the same coefficient
sequence. These identities follow from classical facts and are
independent of the Riemann Hypothesis; the rigorous connection between
these sequences and RH is the Jensen--P\'olya program \cite{Griffin} and
Li's criterion \cite{Li}, within which RH remains open.

\section*{Acknowledgements}

We thank Jacques G\'elinas for a detailed critique of the first version,
which identified the circular use of RH, the divergent Bernoulli series,
and the availability of more accurate coefficient data, and which
prompted the present corrections.

\end{document}